\input amstex
\documentstyle{amsppt}
\pagewidth{5.4in}
\pageheight{7.6in}
\magnification=1200
\TagsOnRight
\NoRunningHeads
\topmatter
\title
\bf Maximum principle and convergence of\\
fundamental solutions for the Ricci flow
\endtitle
\author
Shu-Yu Hsu
\endauthor
\affil
Department of Mathematics\\
National Chung Cheng University\\
168 University Road, Min-Hsiung\\
Chia-Yi 621, Taiwan, R.O.C.\\
e-mail:syhsu\@math.ccu.edu.tw
\endaffil
\date
Nov 8, 2007
\enddate
\address
e-mail address:syhsu\@math.ccu.edu.tw
\endaddress
\abstract
In this paper we will prove a maximum principle for the solutions
of linear parabolic equation on complete non-compact manifolds with 
a time varying metric. We will prove the convergence of the Neumann 
Green function of the conjugate heat equation for the Ricci flow 
in $B_k\times (0,T)$ to the minimal fundamental solution of the 
conjugate heat equation as $k\to\infty$. We will prove the uniqueness of 
the fundamental solution under some exponential decay assumption on 
the fundamental solution. We will also give a detail proof of the 
convergence of the fundamental solutions of the conjugate heat equation 
for a sequence of pointed Ricci flow $(M_k\times (-\alpha,0],x_k,g_k)$ 
to the fundamental solution of the limit manifold as $k\to\infty$ which 
was used without proof by Perelman in his proof of the pseudolocality 
theorem for Ricci flow \cite{P}. 
\endabstract
\endtopmatter
\NoBlackBoxes
\define \pd#1#2{\frac{\partial #1}{\partial #2}}
\define \1{\partial}
\define \2{\overline}
\define \3{\varepsilon}
\define \4{\widetilde}
\define \5{\underline}
\define \ov#1#2{\overset{#1}\to{#2}}
\define \oa#1{\overset{a}\to{#1}}
\document

Maximum principle for the heat equation on complete non-compact manifold
with a fixed metric was proved by P.~Li, L.~Karp \cite{LK} and J.~Wang 
\cite{W} (cf. \cite{CLN}). Maximum principle for parabolic equations on 
complete non-compact manifold with a metric with uniformly bounded Riemannian 
curvature and evolving by the Ricci flow,
$$ 
\frac{\1 }{\1 t}g_{ij}=-2R_{ij}\tag 0.1
$$ 
was proved by W.X.~Shi \cite{S1}, \cite{S2}, \cite{S3} under either
a uniform boundedness condition on the solution or some structural 
conditions on the parabolic equation or positivity assumption on the 
Riemannian curvature operator. 

Let $M$ be a non-compact manifold with a time varying metric
$g(t)=(g_{ij}(t))$, $0\le t<T$, such that for any $0\le t<T$ $(M, g(t))$ is 
a complete non-compact manifold. Let $x_0$ be a fixed point of $M$.
In this paper we will prove the maximum principle for
the subsolution of the linear parabolic equation,
$$
u_t=\Delta u+\vec{a}\cdot\nabla u+bu\quad\text{ in } M\times (0,T) 
\tag 0.2
$$
under the condition
$$
\int_0^T\int_Mu_+^2(x,t)e^{-\lambda r_t(x_0,x)^2}\,dV_t\,dt<\infty
\tag 0.3
$$
for some constant $\lambda>0$, vector field $\vec{a}(\cdot,t)$, $0\le t\le T$,
and function $b(x,t)$ on $M\times [0,T)$ where $u_+=\max (u,0)$ and 
$r_t(x_0,x)$ is the distance between $x_0$ and $x$ with respect to the 
metric $g_{ij}(t)$.
 
In \cite{CTY} A. Chau, L.F.~Tam and C.~Yu proved the existence of minimal 
fundamental solution of the conjugate heat equation of Ricci flow on any
n-dimensional non-compact complete manifold, $n\ge 3$, by approximating it 
by a monotone increasing sequence of Dirichlet Green functions of the 
conjugate heat equation of Ricci flow in bounded domains. In this paper we 
will show that their argument can be modified to work for any $n\ge 2$. We
will prove that the Neumann Green functions of the conjugate heat equation 
of Ricci flow in bounded domains will also converge to the minimal 
fundamental solution of the conjugate heat equation of \cite{CTY} for any
$n\ge 2$. 

We will prove the uniqueness of the fundamental solution of the conjugate 
heat equation under some exponential decay assumption on the fundamental 
solution. We will also give a detail proof of the 
convergence of the fundamental solutions of the conjugate heat equation 
for a sequence of pointed Ricci flow $(M_k\times (-\alpha,0],x_k,g_k)$ 
to the fundamental solution of the limit manifold as $k\to\infty$ which 
was used without proof by Perelman in his proof of the pseudolocality 
theorem for Ricci flow \cite{P}. 

We start will some definitions. Let $x_0\in M$ and let $r(x,y)=r_0(x,y)$, 
$r_t(x)=r_t(x_0,x)$, $r(x)=r(x_0,x)=r_0(x_0,x)$. Let $\nabla^t$ and 
$\Delta^t$ be the 
covariant derivative and Laplacian with respect to the metric $g(t)$. 
When there is no ambiguity, we will drop the superscript and write $\nabla$,
$\Delta$, for $\nabla^t$, $\Delta^t$, respectively. For any $R>0$, 
$y\in M$, let $B_R^t(y)=B_{g(t)}(y,R)$ be the geodesic ball with center 
$y$ and radius $R$ with respect to the metric $g(t)$ and let $B_R=B_R^0(x_0)$.
Let $dV_t$, $dV$, be the volume element with respect to the metric $g(t)$
and $g(0)$ respectively and let $V_x(r)=\text{Vol}_{g(0)}(B_r(x))$.

$$
\text{Section 1}
$$

\proclaim{\bf Theorem 1.1}
Let $M$ be a non-compact manifold with a time varying metric
$g(t)=(g_{ij}(t))$, $0\le t<T$, such that for any $0\le t<T$ $(M, g(t))$ is 
a complete non-compact manifold. Let $\vec{a}(\cdot,t)$, $0\le t<T$,
be a vector field on $M$ which satisfies 
$$
\sup_{M\times [0,T)}|\vec{a}|\le\alpha_1\tag 1.1
$$
and let $b\in L^{\infty}(M\times [0,T))$ such that
$$
\sup_{M\times [0,T)}|b|\le\alpha_2\tag 1.2
$$ 
for some constants $\alpha_1>0$, $\alpha_2>0$. Suppose $g(t)$ satisfies
$$
-\alpha_3g_{ij}\le\frac{\1 g_{ij}}{\1 t}\le\alpha_3g_{ij}
\quad\text{ in }M\times (0,T)\tag 1.3
$$
for some constant $\alpha_3>0$ and $u\in C(M\times [0,T))\cap
C^{2,1}(M\times (0,T))$ is a subsolution of (0.2) satisfying (0.3) for 
some constant $\lambda>0$ and 
$$
u(x,0)\le 0\quad\forall x\in M.\tag 1.4
$$
Then 
$$
u(x,t)\le 0\quad\text{ on }M\times [0,T).\tag 1.5
$$
\endproclaim
\demo{Proof}
We will use a modification of the proof in \cite{EH}, \cite{LK}, \cite{NT} 
and \cite{W} to prove the theorem. Let $x_0\in M$, $r_t(x)=r_t(x_0,x)$, 
$r(x)=r_0(x)$, and 
$$
h(x,t)=-\frac{r(x)^2}{4(2\eta -t)}\quad\forall x\in M,0\le t\le\eta
$$ 
for some constant $0<\eta\le(\log (9/8))/\alpha_3$ to be determined later. 
Then $h$ satisfies
$$
h_t+|\nabla^0h|_{g(0)}^2=0\quad\text{ in }M\times [0,\eta].\tag 1.6
$$
Choose a smooth function $\phi:\Bbb{R}\to\Bbb{R}$, $0\le\phi\le 1$, such 
that $\phi (x)=1$ for all $x\le 0$, $\phi (x)=0$ for all $x\ge 1$ and 
$-2\le\phi'(x)\le 0$ for any $x\in\Bbb{R}$. For any $R\ge 1$, let $\phi_R(x)
=\phi (r(x)-R)$. Then $|\nabla^0\phi_R|_{g(0)}\le 2$ on $M$. Now by (1.3),
$$\align
&\left\{\aligned
&e^{-\alpha_3t}g(0)\le g(t)\le e^{\alpha_3t}g(0)\qquad
\qquad\text{ in }M\times [0,T)\\
&e^{-\alpha_3t}g^{-1}(0)\le g^{-1}(t)\le e^{\alpha_3t}g^{-1}(0)
\quad\text{ in }M\times[0,T)\endaligned\right.\tag 1.7\\
\Rightarrow\quad&e^{-\alpha_3T/2}r(x)\le r_t(x)\le e^{\alpha_3T/2}r(x)
\qquad\quad\forall x\in M,0\le t<T\tag 1.8
\endalign
$$
and
$$\align
&\biggl |\frac{\1}{\1 t}(dV_t)\biggr |\le\frac{n\alpha_3}{2}dV_t
\qquad\qquad\quad\text{ in } M\quad\forall 0\le t<T\\
\Rightarrow\quad&e^{-\frac{n\alpha_3}{2}T}V_s\le dV_t\le 
e^{\frac{n\alpha_3}{2}T}dV_s\quad\text{ in } M\quad\forall 0\le s,t<T.
\tag 1.9
\endalign
$$
Hence $|\nabla^t\phi_R|\le 2e^{\alpha_3T/2}$ on $M\times [0,T)$.
By (1.6) and (1.7),
$$
h_t+e^{-\alpha_3\eta}|\nabla h|^2\le 0\quad\text{ in }M\times [0,\eta].
\tag 1.10
$$
Then by (0.2), (1.9) and (1.10),
$$\align
&\frac{\1}{\1 t}\biggl (\int_M\phi_R^2e^hu_+^2\,dV_t\biggr )\\
=&\int_M\phi_R^2e^hh_tu_+^2\,dV_t+2\int_M\phi_R^2e^hu_+u_t\,dV_t
+\int_M\phi_R^2e^hu_+^2\,\frac{\1}{\1 t}(dV_t)\\
\le&\int_M\phi_R^2e^hh_tu_+^2\,dV_t+2\int_M\phi_R^2e^hu_+\Delta u\,dV_t
+2\int_M\phi_R^2e^hu_+\vec{a}\cdot\nabla u\,dV_t
+2\int_M\phi_R^2e^hbu_+^2\,dV_t\\
&\qquad +\frac{n\alpha_3}{2}\int_M\phi_R^2e^hu_+^2\,dV_t\\
\le&-e^{-\alpha_3\eta}\int_M\phi_R^2e^h|\nabla h|^2u_+^2\,dV_t
-2\int_M\phi_R^2e^h|\nabla u_+|^2\,dV_t-2\int_M\phi_R^2e^hu_+\nabla h\cdot
\nabla u_+\,dV_t\\
&\qquad-4\int_M\phi_Re^hu_+\nabla\phi_R\cdot\nabla u_+\,dV_t
+2\alpha_1\int_M\phi_R^2e^hu_+|\nabla u_+|\,dV_t\\
&\qquad +\biggl (2\alpha_2+\frac{n\alpha_3}{2}\biggr )
\int_M\phi_R^2e^hu_+^2\,dV_t\qquad\qquad\qquad\qquad\qquad
\forall 0\le t<\eta.\tag 1.11
\endalign
$$
Now $\forall 0\le t<\eta$,
$$\align
2\biggl |\int_M\phi_R^2e^hu_+\nabla h\cdot\nabla u_+\,dV_t\biggr |
\le&e^{-\alpha_3\eta}\int_M\phi_R^2e^h|\nabla h|^2u_+^2\,dV_t
+e^{\alpha_3\eta}\int_M\phi_R^2e^h|\nabla u_+|^2\,dV_t\\
\le&e^{-\alpha_3\eta}\int_M\phi_R^2e^h|\nabla h|^2u_+^2\,dV_t
+\frac{9}{8}\int_M\phi_R^2e^h|\nabla u_+|^2\,dV_t,\tag 1.12
\endalign
$$
$$\align
&4\biggl |\int_M\phi_Re^hu_+\nabla\phi_R\cdot\nabla u_+\,dV_t\biggr |\\
\le&\frac{1}{2}\int_M\phi_R^2e^h|\nabla u_+|^2\,dV_t
+8\int_Me^h|\nabla\phi_R|^2u_+^2\,dV_t\\
\le&\frac{1}{2}\int_M\phi_R^2e^h|\nabla u_+|^2\,dV_t
+32e^{\alpha_3T}\int_{B_{R+1}\setminus B_R}e^hu_+^2\,dV_t
\tag 1.13
\endalign
$$
and
$$
2\alpha_1\biggl |\int_M\phi_R^2e^hu_+|\nabla u_+|\,dV_t\biggr |
\le\frac{1}{4}\int_M\phi_R^2e^h|\nabla u_+|^2\,dV_t
+4\alpha_1^2\int_M\phi_R^2e^hu_+^2\,dV_t.\tag 1.14
$$
By (1.11), (1.12), (1.13) and (1.14),
$$\align
&\frac{\1}{\1 t}\biggl (\int_M\phi_R^2e^hu_+^2\,dV_t\biggr )\\
\le&-\frac{1}{8}\int_M\phi_R^2e^h|\nabla u_+|^2\,dV_t
+C_1\int_M\phi_R^2e^hu_+^2\,dV_t+32e^{\alpha_3T}
\int_{B_{R+1}\setminus B_R}e^hu_+^2\,dV_t\\
\Rightarrow\quad&\frac{\1}{\1 t}\biggl (e^{-C_1t}\int_M\phi_R^2e^hu_+^2
\,dV_t\biggr )+\frac{e^{-C_1t}}{8}\int_M\phi_R^2e^h|\nabla u_+|^2\,dV_t\\
\le&32e^{\alpha_3T}\int_{B_{R+1}\setminus B_R}
e^hu_+^2\,dV_t\qquad\qquad\qquad\qquad\qquad\qquad\forall 0\le t<\eta\\
\Rightarrow\quad&e^{-C_1t}\int_M\phi_R^2e^hu_+^2\,dV_t
+\frac{e^{-C_1\eta}}{8}\int_0^t\int_M\phi_R^2e^h|\nabla u_+|^2\,dV_t\,dt\\
\le&32e^{\alpha_3T}\int_0^{\eta}\int_{B_{R+1}\setminus B_R}
e^hu_+^2\,dV_t\,dt\qquad\qquad\qquad\qquad\quad\forall 0\le t<\eta\tag 1.15
\endalign
$$
where $C_1=2\alpha_2+4\alpha_1^2+(n\alpha_3/2)$. By (0.3) and (1.8),
$$
\int_0^{\eta}\int_Mu_+^2(x,t)e^{-\lambda_1 r(x)^2}\,dV_t\,dt<\infty
\tag 1.16
$$
where $\lambda_1=\lambda e^{\alpha_3T}$. We now choose
$\eta=\min (1/(8\lambda_1),(\log (9/8))/\alpha_3)$. Then
$$
h(x,t)\le-\lambda_1r(x)^2\quad\forall x\in M,0\le t<\eta.\tag 1.17
$$
By (1.16) and (1.17),
$$
\int_0^{\eta}\int_Me^hu_+^2(x,t)\,dV_t\,dt<\infty.
\tag 1.18
$$
Letting $R\to\infty$ in (1.15), by (1.18) we get
$$\align
&e^{-C_1t}\int_Me^hu_+^2\,dV_t
+\frac{e^{-C_1\eta}}{8}\int_0^{\eta}\int_Me^h|\nabla u_+|^2\,dV_t\,dt=0
\quad\forall x\in M,0\le t<\min (T,\eta)\\
\Rightarrow\quad&u_+(x,t)=0\qquad\qquad\qquad\qquad\qquad\qquad\qquad\qquad
\qquad\qquad\forall x\in M,0\le t<\min (T,\eta).
\endalign
$$
If $T\le\eta$, we are done. If $T>\eta$, we repeat the above
argument a finite number of times and the theorem follows.
\enddemo

\proclaim{\bf Corollary 1.2}(Lemma 6.2 of \cite{CTY})
Let $(M,g(t))$ be a complete solution of the Ricci flow (0.1) in $(0,T)$ 
with 
$$
|\text{Rm}|\le k_0\quad\text{ on }M\times [0,T)\tag 1.19
$$
for some constant $k_0>0$. Let $u\in C(M\times [0,T))\cap C^{2,1}(M\times 
(0,T))$ satisfy
$$
\Delta u\ge u_t\quad\text{ in }M\times (0,T)
$$
and (0.3), (1.4), for some constant $\lambda>0$. Then $u$
satisfies (1.5). 
\endproclaim

\proclaim{\bf Corollary 1.3}
Let $(M,g(t))$ with $0\le t<T$, $\vec{a}$, $b$ and $\alpha_3$ be as given 
in Theorem 1.1. Suppose (1.19) holds for some constant $k_0>0$.
Let $u\in L^{\infty}(M\times [0,T))\cap C(M\times [0,T))
\cap C^{2,1}(M\times (0,T))$ be a subsolution of (0.2) in $M\times (0,T)$ 
which satisfies (1.4). Then $u$ satisfies (1.5).  
\endproclaim
\demo{Proof}
By the proof of Theorem 1.1 $u$ satisfies (1.8). By the volume 
comparison theorem \cite{C}, there exist 
constants $c_1>0$, $c_2>0$ such that for any $\lambda>0$, $R>0$,
$$
\int_0^T\int_{B_R}u_+^2(x,t)e^{-\lambda r_t(x_0,x)^2}\,dV_t\,dt
\le c_1\|u\|_{\infty}^2\int_0^T\int_0^{\infty}
e^{c_2r-\lambda e^{-\alpha_3T}r^2}\,dr\,dt<\infty.
$$
Letting $R\to\infty$, we get (0.3). Thus the corollary follows from 
Theorem 1.1.
\enddemo

\proclaim{\bf Theorem 1.4}
Let $M$ be a n-dimensional non-compact manifold, $n\ge 2$, such that 
$(M,g(t))$ is a complete solution of the backward Ricci flow 
$$ 
\frac{\1 }{\1 t}g_{ij}=2R_{ij}\tag 1.20
$$  
in $[0,T]$ which satisfies (1.19) for some constant $k_0>0$. Let 
$\Cal{Z}(x,t;y,s)$, $x,y\in M$, $0\le s<t\le T$, be the minimal  
fundamental solution of the forward conjugate heat equation 
in $M\times (s,T]$. That is $\Cal{Z}(\cdot,\cdot;y,s)$ satisfies 
(cf. \cite{CTY})
$$
\1_tu=\Delta u-Ru\quad\text{ in }M\times (s,T]\tag 1.21
$$
with
$$
\lim_{t\searrow s}\Cal{Z}(x,t;y,s)=\delta_y(x).\tag 1.22
$$
For any $k\in\Bbb{Z}^+$, let $\Cal{Z}_k=\Cal{Z}_k(x,t;y,s)$, 
$x,y\in M$, $0\le s<t\le T$, be the Neumann Green function of  
the forward conjugate heat equation which satisfies 
$$\left\{\aligned
&\1_t\Cal{Z}_k=\Delta\Cal{Z}_k-R\Cal{Z}_k\quad\text{ in }B_k\times (s,T]\\
&\frac{\1\Cal{Z}_k}{\1\nu}=0\qquad\qquad\qquad\text{on }\1 B_k\times (s,T]\\
&\lim_{t\searrow s}\Cal{Z}_k(x,t;y,s)=\delta_y(x)\endaligned\right.
\tag 1.23
$$
where $\1/\1\nu$ is the derviative with respect to the unit outward normal on 
$\1 B_k\times (s,T]$ and $B_k=B_k(x_0)$ for some fix point $x_0\in M$.
Then $\Cal{Z}_k(x,t;y,s)$ converges uniformly on every compact
subset of $M\times (s,T]$ to $\Cal{Z}(x,t;y,s)$ as $k\to\infty$.
\endproclaim
\demo{Proof}
By (1.20) and (1.23),
$$\align
&\frac{\1}{\1 t}\int_{B_k}\Cal{Z}_k(x,t;y,s)\,dV_t(x)
=\int_{B_k}\Delta\Cal{Z}_k(x,t;y,s)\,dV_t(x)=0\quad\forall 0\le s<t\le T,
k\in\Bbb{Z}^+\\
\Rightarrow\quad&\int_{B_k}\Cal{Z}_k(x,t;y,s)\,dV_t(x)
=\lim_{t'\searrow s}\int_{B_k}\Cal{Z}_k(x,t';y,s)\,dV_t(x)=1\quad\forall 
0\le s<t\le T, k\in\Bbb{Z}^+\tag 1.24
\endalign
$$
Let $R>1$ and fix $(y,s)\in M\times [0,T)$. By (1.19) and (1.20), 
(1.7) and (1.9) holds with $\alpha_3=2(n-1)k_0$.
Let $G_k(x,t;y,s)$ be the Dirichlet Green function of (1.21) in 
$M\times (0,T)$. We now divide the proof into two cases.

\noindent{\bf Case 1:} $n\ge 3$
 
By (1.7), (1.9), (1.24), and Lemma 3.1 of \cite{CTY}, there exists a 
constant $C_1>0$ such that for any $s<t_1<T$, $k\ge 3R$,
$$\align
\Cal{Z}_k(x,t;y,s)\le&\frac{C_1}{r_1^2V_x(r_1)}
\int_{t-4r_1^2}^{t}\int_{B_{2r_1}(x)}\Cal{Z}_k(z,t;y,s)\,dV_0(z)\, dt\\
\le&\frac{C_1'}{r_1^2\min_{z\in B_R}V_z(r_1)}
\int_{t-4r_1^2}^{t}\int_M\Cal{Z}_k(z,t;y,s)\,dV_t(z)\, dt\\
\le&\frac{4C_1'}{\min_{z\in B_R}V_z(r_1)}
\quad\forall x\in\2{B}_R, t_1\le t\le T
\endalign
$$
where $r_1=\min (1/2,\sqrt{t_1-s}/4)$. Hence the sequence 
$\{\Cal{Z}_k(\cdot,\cdot;y,s)\}$ are uniformly bounded on 
$\2{B}_R\times [t_1,T]$ for any $s<t_1<T$, $k\ge 3R$. By (1.23) and the 
parabolic Schauder estimates \cite{LSU} the sequence
$\{\Cal{Z}_k\}$ is uniformly bounded in $C^{2,\beta}(\2{B}_R\times [t_1,T])$
for some $\beta\in (0,1)$ for any $s<t_1<T$, $k\ge 3R$. 

Let $\{k_i\}_{i=1}^{\infty}$ be a sequence of positive integers such that 
$k_i\to\infty$ as $i\to\infty$. Then by the Ascoli
Theorem and a diagonalization argument the sequence 
$\{\Cal{Z}_{k_i}\}_{i=1}^{\infty}$ has
a subsequence which we may assume without loss of generality to be the 
sequence itself which converges uniformly on every compact subset of
$M\times (s,T]$ to some solution $\4{\Cal{Z}}(\cdot,\cdot;y,s)$ of (1.21) in 
$M\times (s,T]$ as $k\to\infty$.

By the construction of $\Cal{Z}(x,t;y,s)$ in \cite{CTY} $G_k$ increases 
monotonically to $\Cal{Z}$ as $k\to\infty$. Let $0\le s<t\le T$. By the 
maximum principle,
$$\align
&G_{k_i}(x,t;y,s)\le\Cal{Z}_{k_i}(x,t;y,s)\quad\forall x,y\in B_{k_i},
0\le s<t\le T,i\in\Bbb{Z}^+\tag 1.25\\
\Rightarrow\quad&\Cal{Z}(x,t;y,s)\le\4{\Cal{Z}}(x,t;y,s)
\quad\forall x,y\in M,0\le s<t\le T\quad\text{ as }i\to\infty.\tag 1.26
\endalign
$$
By (1.24) and (1.25), $\forall 0\le s<t\le T$,
$$\align
&\int_{B_R}G_{k_i}(x,t;y,s)\,dV_t(x)
\le\int_{B_R}\Cal{Z}_{k_i}(x,t;y,s)\,dV_t(x)\le 1\quad\forall k_i>R>1\\
\Rightarrow\quad&\int_{B_R}\Cal{Z}(x,t;y,s)\,dV_t(x)
\le\int_{B_R}\4{\Cal{Z}}(x,t;y,s)\,dV_t(x)\le 1\quad\forall R>1
\text{ as }i\to\infty\\
\Rightarrow\quad&\int_M\Cal{Z}(x,t;y,s)\,dV_t(x)
\le\int_M\4{\Cal{Z}}(x,t;y,s)\,dV_t(x)\le 1\quad\text{ as }R\to\infty
\tag 1.27
\endalign
$$
By Lemma 5.1 of \cite{CTY},
$$
\int_M\Cal{Z}(x,t;y,s)\,dV_t(x)=1\quad\forall y\in M,0\le s<t\le T.\tag 1.28
$$
By (1.26), (1.27), and (1.28),
$$\align
&\int_M\4{\Cal{Z}}(x,t;y,s)\,dV_t(x)=\int_M\Cal{Z}(x,t;y,s)\,dV_t(x)=1
\quad\forall y\in M,0\le s<t\le T\\
\Rightarrow\quad&\4{\Cal{Z}}(x,t;y,s)\equiv\Cal{Z}(x,t;y,s)\qquad\qquad
\qquad\qquad\qquad\qquad\forall x,y\in M,0\le s<t\le T.
\endalign
$$
Since the sequence $\{k_i\}$ is arbitrary, $\Cal{Z}_k(x,t,y,s)$ converges to
$\Cal{Z}(x,t,y,s)$ uniformly on every compact subset of $M\times (s,T]$ as
$k\to\infty$.

\noindent{\bf Case 2:} $n=2$

Let $h_0$ be the standard metric on the $2-$sphere $S^2$ with constant
scalar curvature $1$. Then $(S^2,h(t))$, $h=(h_{\alpha\beta})$, 
$0\le t\le T$, with $h(t)=(1+t)h_0$ is the solution of the backward Ricci 
flow on $S^2$ (P.65 of \cite{MT}). Consider the
manifold $\4{M}=M\times S^2$ with metric $\4{g}$ given by
$\4{g}_{ij}=g_{ij}$, $\4{g}_{\alpha\beta}=h_{\alpha\beta}$, 
$\4{g}_{i\alpha}=\4{g}_{\alpha i}=0$. Then $(\4{M},\4{g})$ satisfies
the backward Ricci flow on $[0,T]$ with uniformly bounded Riemannian
curvatures on $[0,T]$.     

Hence as before there exist constants $C_2>0$, $C_3>0$, such that
$$\left\{\aligned
&\frac{1}{C_2}\4{g}(0)\le\4{g}(t)\le C_2\4{g}(0)\quad\text{ in }\4{M}
\times [0,T]\\
&\frac{1}{C_3}d\4{V}\le d\4{V}_t\le C_3d\4{V}\quad\text{ in }\4{M}
\times [0,T]\endaligned\right.
\tag 1.29
$$
where $d\4{V}_t$ is the volume element of $\4{M}$ with respect to the metric
$\4{g}(t)$ and $d\4{V}=d\4{V}_0$ . For any $x,y\in M$, $x',y'\in S^2$, let
$$
\4{\Cal{Z}}_k(x,x',t;y,y',s)=\Cal{Z}_k(x,t;y,s)\quad\forall 0\le s<t\le T.
$$
Since $\Cal{Z}_k$ satisfies (1.21),
$$
\1_t\4{\Cal{Z}}_k=\Delta_{\4{g}(t)}\4{\Cal{Z}}_k-R\4{\Cal{Z}}_k
\quad\text{ in }B_k\times S^2\times (s,T).\tag 1.30
$$
Then by Lemma 3.1 of \cite{CTY}, (1.24), (1.29) and (1.30), there exists 
a constant $C_4>0$ such that for any $0\le s<t_1\le T$, $k\ge 3R$,
$$\align
\4{\Cal{Z}}_k(x,x',t;y,y',s)\le&\frac{C_4}{r_1^2\4{V}_{(x,x')}(r_1)}
\int_{t-4r_1^2}^{t}\int_{\4{B}_{2r_1}(x,x')}\4{\Cal{Z}}_k(z,t;y,y',s)
\,d\4{V}(z)\, dt\\
\le&\frac{C_4'}{r_1^2\4{V}_{(x,x')}(r_1)}
\int_{t-4r_1^2}^{t}\int_{B_k}\Cal{Z}_k(w,t;y,s)\,dV_t(w)\, dt\\
\le&\frac{4C_4'}{\4{V}_{(x,x')}(r_1)}
\quad\forall x,y\in\2{B}_R,x',y'\in S^2,t_1\le t\le T\\
\Rightarrow\qquad\quad\Cal{Z}_k(x,t;y,s)
\le&\frac{4C_4'}{\min\Sb w\in\2{B}_R\\x'\in S^2\endSb
\4{V}_{(w,x')}(r_1)}\quad\forall x\in B_R,t_1\le t\le T\tag 1.31
\endalign
$$
where $r_1=\min (1/2,\sqrt{t_1-s}/4)$, $\4{B}_{2r_1}(x,x')$ is the geodesic 
ball of radius $2r_1$ and center $(x,x')$ in $\4{M}$ with respect to the 
metric $\4{g}(0)$ and $\4{V}_{(x,x')}(r_1)$ is the 
volume of $\4{B}_{r_1}(x,x')$ with respect to the metric $\4{g}(0)$. Hence the 
sequence $\{\Cal{Z}_k(\cdot,\cdot;y,s)\}_{k=1}^{\infty}$ are uniformly 
bounded on $\2{B}_R\times [t_1,T]$ for any $s<t_1\le T$, $k\ge 3R$. 

By a similar argument the sequence $\{G_k(\cdot,\cdot;y,s)\}_{k=1}^{\infty}$ 
are uniformly bounded on $\2{B}_R\times [t_1,T]$ for any 
$s<t_1\le T$, $k\ge 3R$. Then by the same argument as in \cite{CTY},
$G_k$ increases monotonically to $\Cal{Z}$ as $k\to\infty$ and
$\Cal{Z}(x,t;y,s)$ satisfies (1.28). By (1.28), (1.31), and an argument 
similar to case 1, $\Cal{Z}_k(x,t;y,s)$ converges to $\Cal{Z}(x,t;y,s)$ 
uniformly on every compact subset of $M\times (s,T]$ as $k\to\infty$ and 
the theorem follows.
\enddemo

\proclaim{\bf Corollary 1.5}
Let $(M,g(t))$, $0\le t\le T$, and $\Cal{Z}(x,t;y,s)$, $x,y\in M$,
$0\le s<t\le T$, be as in Theorem 1.4. Suppose $\4{\Cal{Z}}(x,t;y,s)$ is a 
fundamental solution of the forward conjugate heat equation which satisfies 
(1.21), (1.22) and
$$
\forall y\in M,\max_{s\le t\le T}\int_{B_R}\4{\Cal{Z}}(x,t;y,s)\,dV_t(x)
\le o(R)\quad\forall 0\le s<T\quad\text{ as }R\to\infty,\tag 1.32
$$
then 
$$
\4{\Cal{Z}}(x,t;y,s)\equiv\Cal{Z}(x,t;y,s)\quad\forall x,y\in M,
0\le s<t\le T.\tag 1.33
$$
\endproclaim
\demo{Proof}
By (1.32) and an argument similar to the proof of Lemma 5.1 of \cite{CTY},
$$
\int_M\4{\Cal{Z}}(x,t;y,s)\,dV_t(x)=1\quad\forall y\in M,
0\le s<t\le T.\tag 1.34
$$
Let $G_k$ be as in the proof of Theorem 1.4. By the maximum principle,
$$
G_k(x,t;y,s)\le\4{\Cal{Z}}(x,t;y,s)\quad\forall x,y\in B_k,0\le s<t\le T,
k\in\Bbb{Z}^+.\tag 1.35
$$
Since $\Cal{Z}$ satisfies (1.28), by (1.34), (1.35) and an argument 
similar to the proof of Theorem 1.4 the corollary follows.
\enddemo

\proclaim{\bf Theorem 1.6}
Let $(M,g(t))$, $0\le t\le T$, and $\Cal{Z}(x,t;y,s)$, $x,y\in M$,
$0\le s<t\le T$, be as in Theorem 1.4. Let $\4{\Cal{Z}}(x,t;y,s)$ be 
a fundamental solution of the forward conjugate heat equation which 
satisfies (1.21) and (1.22). Then (1.33) holds if and only if there 
exist constants $C>0$ and $D>0$ such that
$$
\4{\Cal{Z}}(x,t;y,s)\le\frac{C}{V_y(\sqrt{t-s})}
e^{-\frac{r^2(x,y)}{D(t-s)}}\quad\forall 0\le s<t\le T.\tag 1.36
$$
\endproclaim
\demo{Proof}
The case (1.33) implies (1.36) was proved in \cite{CTY}. Hence we only 
need to show that (1.36) implies (1.33). Suppose there exist constants 
$C>0$ and $D>0$ such that (1.36) holds. By the proof of Theorem 1.4
(1.9) holds for some constant $\alpha_3>0$. Then by (1.9) and (1.36),
$$\align
\int_M\4{\Cal{Z}}(x,t;y,s)\,dV_t(x)\le&\frac{C}{V_y(\sqrt{t-s})}
\int_Me^{-\frac{r^2(x,y)}{D(t-s)}}\,dV_t(x)\\
\le&\frac{C}{V_y(\sqrt{t-s})}\int_0^{\infty}e^{-\frac{r^2}{D(t-s)}}\,dV_y(r)\\
=&C\int_0^{\infty}\frac{V_y(r)}{V_y(\sqrt{t-s})}e^{-\frac{r^2}{D(t-s)}}
\,d\biggl (\frac{r^2}{D(t-s)}\biggr).\tag 1.37
\endalign
$$
Let $V_{k_0}(r)$ be the volume of the geodesic ball of radius $r$ 
in the space form with constant sectional curvature $-k_0$. 
Let $\delta=\sqrt{t-s}$ and $a=(r/\sqrt{t-s})+1$. By (1.19) and the volume 
comparison theorem \cite{C},
$$\align
\frac{V_y(r)}{V_y(\sqrt{t-s})}\le&\frac{V_y(r+\sqrt{t-s})}{V_y(\sqrt{t-s})}
\le\frac{V_{k_0}(r+\sqrt{t-s})}{V_{k_0}(\sqrt{t-s})}
=\frac{V_{k_0}(a\delta)}{V_{k_0}(\delta)}
=a^n\frac{V_{a^2k_0}(\delta)}{V_{k_0}(\delta)}\\
\le&a^n\frac{V_{a^2k_0}(\sqrt{T})}{V_{k_0}(\sqrt{T})}
=\frac{V_{k_0}(a\sqrt{T})}{V_{k_0}(\sqrt{T})}.\tag 1.38
\endalign
$$
Now
$$\align
V_{k_0}(a\sqrt{T})=&\int_0^{a\sqrt{T}}\biggl (\frac{1}{\sqrt{k_0}}\sinh 
(\sqrt{k_0}\rho)\biggr )^{n-1}\,d\rho\le Ce^{(n-1)\sqrt{{k_0}T}a}\\
=&Ce^{(n-1)\sqrt{{k_0}T}(\frac{r}{\sqrt{t-s}}+1)}\\
\le&Ce^{C'\frac{r}{\sqrt{t-s}}}.\tag 1.39
\endalign
$$
By (1.37), (1.38), and (1.39),
$$\align
\int_M\4{\Cal{Z}}(x,t;y,s)\,dV_t(x)
\le&C\int_0^{\infty}e^{-\frac{r^2}{D(t-s)}+C'\frac{r}{\sqrt{t-s}}}
\,d\biggl (\frac{r^2}{D(t-s)}\biggr)=C_T<\infty\quad\forall
0\le s<t\le T
\endalign
$$
for some constant $C_T>0$ depending on $k_0$ and $T$. Hence by Corollary 1.5 
(1.33) holds.
\enddemo

$$
\text{Section 2}
$$

In this section we will give a detail proof of the convergence of fundamental
solutions of conjugate heat equation for Ricci flow which was used without 
proof by Perelman in his proof of the pseudolocality theorem for Ricci flow 
\cite{P}. 

\proclaim{\bf Theorem 2.1}
Let $\alpha>0$ and let $(M_k\times (-\alpha,0], x_k, g_k(t))$ 
be a sequence of pointed Ricci flow (0.1) where each $M_k$ is either a 
closed manifold or a non-compact manifold with bounded curvature
such that $(M_k,g_k(t))$ is complete for each $-\alpha<t\le 0$. Suppose 
$$
|\text{Rm}_k|(x,t)\le C_1\quad\forall x\in B_k(x_k,A_k),
-\alpha<t\le 0, k\in\Bbb{Z}^+\tag 2.1
$$
for some constant $C_1>0$ and sequence $\{A_k\}$, $A_k\to\infty$ as
$k\to\infty$, 
and 
$$(M_k\times (-\alpha,0], x_k, g_k(t))
$$  
converges in the $C^{\infty}$-sense to some pointed Ricci flow 
$(M\times (-\alpha,0],x_{\infty},g_{\infty})$ as $k\to\infty$ where
$B_k(x_k,A_k)=B_{g_k(0)}(x_k,A_k)$.
That is there exists an exhausting sequence $U_1\subset U_2\subset\cdots
\subset M$ of open sets each containing $x_{\infty}$ and each with 
compact closure in $M$ and diffeomeomorphisms $\Phi_k$ of $U_k$ to open sets
$V_k$ of $M_k$ such that $\Phi_k(x_{\infty})=x_k\quad\forall 
k\in\Bbb{Z}^+$ and the pull-back metric $\Phi_k^{\ast}(g_k)$ 
converges uniformly to $g_{\infty}$ on every compact subset of $M\times 
(-\alpha,0]$ as $k\to\infty$. 

If $u_k$ satisfies the conjugate heat equation,
$$
u_t+\Delta_ku-R_{g_k(t)}u=0
\tag 2.2
$$
in $M_k\times (-\alpha,0)$ with
$$
\lim_{t\nearrow 0}u(x,t)=\delta_{x_k}
$$
where $\Delta_k=\Delta_{g_k(t)}$,
then $\Phi_k^{\ast}(u_k)$ will converge uniformly on every compact subset 
of $M\times (-\alpha,0)$ to the minimal fundamental solution $u$ of the 
conjugate heat equation 
$$
u_t+\Delta_{g_{\infty}(t)}u-R_{g_{\infty}(t)}u=0\tag 2.3
$$
of $(M,g_{\infty})$ in $M\times (-\alpha,0)$ with
$$
\lim_{t\nearrow 0}u(x,t)=\delta_{x_{\infty}}\tag 2.4
$$  
as $k\to\infty$.
\endproclaim
\demo{Proof}
For simplicity we will write $B_r(x)$, $V_k$, $V$, $dV_k^t$, $dV_t$, for 
$B_{g_{\infty}(0)}(x,r)$, $\text{Vol}_{g_k(0)}$, 
$\text{Vol}_{g_{\infty}(0)}$, $dV_{g_k(t)}$, and $dV_{g_{\infty}(t)}$ 
respectively. We also let $dV_k=dV_k^0$. Note that $u_k>0$ in $M_k\times 
(-\alpha,0)$ and (\cite{CTY})
$$
\int_{M_k}u_k(y,t)\,dV_k^t(y)=1\quad\forall -\alpha<t\le 0, 
k\in\Bbb{Z}^+.\tag 2.5
$$
Since $g_k$ satisfies (0.1) in $M_k\times (-\alpha,0]$, by (2.1) there 
exists a constant $C_2>1$ such that
$$\left\{\aligned
&\frac{1}{C_2}g_k(x,s)\le g_k(x,t)\le C_2g_k(x,s)\quad\forall 
x\in B_k(x_k,A_k),-\alpha<s,t\le 0,k\in\Bbb{Z}^+\\
&\frac{1}{C_2}dV_k^s(x)\le dV_k^t(x)\le C_2dV_k^s(x)\quad\forall 
x\in B_k(x_k,A_k),-\alpha<s,t\le 0,k\in\Bbb{Z}^+.
\endaligned\right.\tag 2.6
$$
Let $r_k(x,y,t)$ be the geodesic distance between $x,y\in M_k$ with respect
to the metric $g_k(t)$ and let $r_k(x,y)=r_k(x,y,0)$. Let
$r(x,y)$ be the geodesic distance between $x,y\in M$ with respect
to the metric $g_{\infty}(0)$.
Let $R>1$, $-\alpha<t_1<t_2<0$, 
and let $r_1=\min (1/2,\sqrt{-t_2}/4)$. We choose $k_1'\in\Bbb{Z}^+$ such that
$\2{B_{R+2}(x_{\infty})}\subset U_{k_1}$ and $A_k\ge 6\sqrt{C_2}R$ 
for all $k\ge k_1'$. Then by (2.6),
$$
\frac{1}{\sqrt{C_2}}r_k(x,y,s)\le r_k(x,y,t)\le\sqrt{C_2}r_k(x,y,s)
\quad\forall x,y\in B_k(x_k,2R),-\alpha<s,t\le 0,k\ge k_1'.
\tag 2.7
$$
Since $\2{B_R(x_{\infty})}\times [t_1,t_2]$ is compact, 
there exist 
$$
z_1, z_2,\dots,z_m\in \2{B_R(x_{\infty})},s_1,s_2, \dots,s_m\in [t_1,t_2]
$$ 
such that 
$$
\2{B_R(x_{\infty})}\times [t_1,t_2]\subset
\cup_{j=1}^mB_{\frac{r_1}{2}}(z_j)\times [s_j,s_j+r_1^2).\tag 2.8
$$  
Let $z_j^k=\Phi_k(z_j)$. Since $\Phi_k^{\ast}(g_k)$ converges 
uniformly to $g_{\infty}$ on $\2{B_{2R}(x_{\infty})}\times [t_1,0]$ 
as $k\to\infty$, there exists $k_2'\ge k_1'$ such that for any $k\ge k_2'$, 
$j=1,\dots,m$, 
$$\left\{\aligned
&\Phi_k(\2{B_R(x_{\infty})})\subset\2{B_k(x_k,R+(1/2))}\\
&\Phi_k(\2{B_{r_1/\sqrt{2C_2}}(z_j)})
\subset\2{B_k(z_j^k,r_1/\sqrt{C_2})}\\
&\Phi_k(\2{B_{\frac{r_1}{2}}(z_j)})\subset\2{B_k(z_j^k,r_1)}.
\endaligned\right.\tag 2.9
$$
By (2.9),
$$
z_j^k\in\2{B_k(x_k,R+(1/2))}\quad\forall j=1,\dots,m,k\ge k_2'\tag 2.10
$$
By (2.7) and (2.10),
$$
B_k^s(z_j^k,C_2^{-\frac{1}{2}}r)\subset B_k^t(z_j^k,r)\quad\forall
t_1\le s,t\le 0,k\ge k_2', j=1,\dots,m,0<r\le r_1.\tag 2.11
$$
Let $\{k_i\}_{i=1}^{\infty}\subset\Bbb{Z}^+$ be a sequence such that 
$k_i\to\infty$ as $i\to\infty$. We now divide the proof into two cases.

\noindent{\bf Case 1:} $n\ge 3$

By (2.1), (2.2), (2.5), (2.6), (2.7), (2.9), (2.10), (2.11) and an argument 
similar to the proof of Lemma 3.1 of \cite{CTY} and Theorem 3.1 of \cite{KZ} 
and there exists a constant $C_3>0$ such that for any 
$x\in\2{B_k(z_j^k,r_1)}$, $s_j\le t\le s_j+r_1^2$, $k\ge k_2'$, $j=1,\dots,m$,
$$\align
u_k(x,t)\le&\frac{C_3}{r_1^2V_k(B_k(z_j^k,C_2^{-\frac{1}{2}}r_1))}
\int_{s_j}^{s_j+4r_1^2}\int_{B_k(z_j^k,2r_1)}u_k(y,t)\,dV_k(y)\,dt\\
\le&\frac{C_2C_3}{r_1^2V_k(B_k(z_j^k,C_2^{-\frac{1}{2}}r_1))}
\int_{s_j}^{s_j+4r_1^2}\int_{B_k(z_j^k,2r_1)}u_k(x,t)\,dV_k^t(y)\,dt\\
\le&\frac{4C_2C_3}{V_{\Phi_k^{\ast}(g_k(0))}(B_{r_1/\sqrt{2C_2}}(z_j))}.
\tag 2.12
\endalign
$$
Since $\Phi_k^{\ast}(g_k)$ converges uniformly to $g_{\infty}$ on 
$\2{B_{2R}(x_{\infty})}\times [t_1,0]$ as $k\to\infty$, there exist
$k_3'\ge k_2'$ and a constant $C_4>0$ such that
$$
V_{\Phi_k^{\ast}(g_k(0))}(B_{r_1/\sqrt{2C_2}}(z_j))
\ge C_4\quad\forall k\ge k_3',j=1,\dots,m.\tag 2.13
$$
By (2.8), (2.9), (2.12) and (2.13),
$$\align
\Phi_k^{\ast}(u_k)(y,t)\le&\frac{4C_2C_3}{C_4}\quad\forall
y\in\2{B_R(x_{\infty})}, t_1\le t\le t_2,k\ge k_3'.\tag 2.14
\endalign
$$
Hence the sequence $\{\Phi_k^{\ast}(u_k)\}_{k=1}^{\infty}$ are uniformly 
bounded on $\2{B_R(x_{\infty})}\times [t_1,t_2]$ for any $-\alpha
<t_1<t_2<0$.
Since $\Phi_k^{\ast}(u_k)$ satisfies the conjugate heat equation on 
$\2{B_R(x_{\infty})}\times (-\alpha,0)$, by (2.2) and the injectivity 
radius estimates of \cite{CLY} and the uniform convergence of $\Phi_k^{\ast}
(g_k)$ to $g_{\infty}$ on every compact subset of $M\times (-\alpha,0]$
as $k\to\infty$, one can apply the parabolic Schauder estimates 
of \cite{LSU} to conclude that for any $R>1$ and $-\alpha<t_1<t_2<0$ 
$\{\Phi_k^{\ast}(u_k)\}_{k=1}^{\infty}$ are uniformly bounded in $C^{2,\beta}
(\2{B_R(x_{\infty})}\times [t_1,t_2])$ for some $\beta\in (0,1)$. 

\noindent{\bf Case 2:} $n=2$

By considering $\4{M}=M\times S^2$ and using an argument similar to the 
proof of Theorem 1.4 and case 1 one can also conclude that
when $n=2$, for any $R>1$ and $-\alpha<t_1<t_2<0$ 
$\{\Phi_k^{\ast}(u_k)\}_{k=1}^{\infty}$ are uniformly 
bounded in $C^{2,\beta}(\2{B_R(x_{\infty})}\times [t_1,t_2])$ for some 
$\beta\in (0,1)$. 

Hence by case 1 and case 2, the Ascoli Theorem and a diagonalization 
argument $\{\Phi_{k_i}^{\ast}(u_{k_i})\}_{i=1}^{\infty}$ has a subsequence 
which we may assume without loss of generality to be the sequence  
$\{\Phi_{k_i}^{\ast}(u_{k_i})\}_{i=1}^{\infty}$  itself
that converges uniformly on every compact subset of $M\times (-\alpha,0)$
to a solution $u$ of the conjugate heat equation of $(M,g_{\infty})$ in 
$M\times (-\alpha,0)$ as $i\to\infty$.

By (2.1), (2.6), (2.7), and an argument similar to the proof of Theorem 5.1 
and (5.2) of \cite{CTY}, there exist constants $C>0$ and $D>0$ such
that
$$\align
&u_k(x,t)\le\frac{C}{V_k(B_k(x_k,\sqrt{t}))}e^{-\frac{r_k(x,x_k)^2}{D|t|}}
\quad\forall x\in B_k(x_k,R),-\alpha<t<0,k\ge k_3'\\
\Rightarrow\quad&u(x,t)\le\frac{C}{V(B_{\sqrt{t}}(x_{\infty}))}
e^{-\frac{r(x,x_{\infty})^2}{D|t|}}
\qquad\forall x\in M,-\alpha<t<0\quad\text{ as }k=k_i\to\infty.
\tag 2.15
\endalign
$$
Let $\psi\in C_0^{\infty}(M)$. Then  supp$\,\psi\subset  
B_{R_1}(x_{\infty})$ for some constant $R_1>0$. Choose $k_4'\ge k_3'$ 
such that $\2{B_{R_1}(x_{\infty})}\subset U_k$ for all $k\ge k_4'$. Let 
$$
\psi_k(x)=\Phi_k^{\ast}(\psi)(x)=\left\{\aligned
&\psi(\Phi_k^{-1}(x))\quad\text{ if }x\in V_k\\
&0\qquad\qquad\quad\text{ if }x\not\in V_k.\endaligned\right.
$$
Then $\psi_k\in C_0^{\infty}(M_k)$ for all $k\ge k_4'$ and $\psi_k(x_k)
=\psi(\Phi_k^{-1}(x_k))=\psi(x_{\infty})$. Let $t_1\in (-\alpha,0)$.
Then by (0.1), (2.2) and (2.5), $\forall t_1\le t<0$, $k\ge k_4'$,
$$\align
\biggl |\int_{M}\Phi_k^{\ast}(u_k)\psi\,dV_{\Phi_k^{\ast}(g_k(t))}
-\psi(x_{\infty})\biggr |
=&\biggl |\int_{M_k}u_k\psi_k\,dV_k^t-\psi_k(x_k)\biggr |\\
=&\biggl |\int_0^t\int_{M_k}\psi_k\biggl (\frac{\1 u_k}{\1 t}
-R_{g_k(t)}u_k\biggr )\,dV_k^t\,dt\biggr |\\
=&\biggl |\int_0^t\int_{M_k}\psi_k\Delta_ku_k\,dV_k^t\,dt\biggr |\\
=&\biggl |\int_0^t\int_{M_k}u_k\Delta_k\psi_k\,dV_k^t\,dt\biggr |\\
\le&\max|\Delta_k\psi_k|\biggl (\int_{M_k}u_k\,dV_k^t\biggr )|t|\\
\le&|t|\max_{\2{B_{R_1}(x_{\infty})}\times [t_1,0]}
|\Delta_{\Phi_k^{\ast}(g_k(t))}\psi|.\tag 2.16
\endalign
$$
Letting $k\to\infty$ in (2.16),
$$\align
&\biggl |\int_{M}u\psi\,dV(t)-\psi(x_{\infty})\biggr |
\le |t|\max_{\2{B_{R_1}(x_{\infty})}\times [t_1,0]}
|\Delta_{g_{\infty}(t)}\psi|\quad\forall t_1\le t<0\\
\Rightarrow\quad&\lim_{t\nearrow 0}\int_{M}u\psi\,dV(t)=\psi(x_{\infty})
\quad\text{ as }t\to 0.
\endalign
$$
Hence $u$ satisfies (2.3) in $M\times (-\alpha,0)$ and (2.4) holds. By (2.15) 
and Theorem 1.6 $u$ is the unique minimal fundamental solution of the 
conjugate heat equation (2.3) in $M\times (-\alpha,0)$ which satisfies (2.4). 
Since the sequence $\{k_i\}_{i=1}^{\infty}$ is arbitrary, 
$\Phi_k^{\ast}(u_k)$ converges uniformly on every compact subset 
of $M\times (-\alpha,0)$ to the minimal fundamental solution of the 
conjugate heat equation of $(M,g_{\infty})$ in $M\times (-\alpha,0)$
which satisfies (2.4) as $k\to\infty$ and the theorem follows.
\enddemo

$$
\text{Acknowledgements}
$$

I woulod like to thank P.~Li for sending me his unpublished paper 
with L.~Karp \cite{LK} and J.~Wang for telling me the maximum principle
results of \cite{W}.


\Refs

\ref
\key CTY\by\ A.~Chau, L.F.~Tam and C.~Yu\paper Pseudolocality 
for the Ricci flow and applications,\ \ \linebreak
http://arxiv.org/abs/math/0701153\endref 

\ref
\key C\by I.~Chavel\book Riemannian geometry:A modern introduction
\publ Cambridge University Press\publaddr Cambridge, United Kingdom
\yr 1995\endref

\ref
\key CLY\by\ \ \ S.Y.~Cheng, P.~Li and S.T.~Yau\paper On the upper estimate 
of the heat kernel of a complete Riemannian manifold\jour Amer. J. Math
\vol 103(5)\pages 1021--1063\yr 1981\endref

\ref
\key CLN\by \ \ \ B.~Chow and P.~Lu and Lei Ni\book Hamilton's Ricci flow,
Graduate Studies in Mathematics, vol. 77\publ Amer. Math. Soc.
\publaddr Providence, R.I., U.S.A.\yr 2006\endref

\ref
\key EH\by\ K.~Ecker and G.~Huisken\paper Interior estimates for 
hypersurfaces moving by mean curvature\jour Invent. Math.\vol 105\yr 1991
\pages 547--569\endref

\ref
\key H1\by R.S.~Hamilton\paper Three-manifolds with positive Ricci curvature
\jour J. Differential Geom.\vol 17(2)\yr 1982\pages 255--306\endref

\ref
\key H2\by R.S.~Hamilton\paper Four-manifolds with positive curvature
operator\jour J. Differential Geom.\vol 24(2)\yr 1986\pages 153--179\endref

\ref 
\key H3\by R.S.~Hamilton\paper The formation of singularities in the Ricci 
flow\jour Surveys in differential geometry, Vol. II (Cambridge, MA, 1993),
7--136, International Press, Cambridge, MA, 1995\endref

\ref
\key H4\by R.S.~Hamilton\paper A compactness property for solutions of the 
Ricci flow\jour Amer. J. Math.\vol 117(3)\yr 1995\pages 545--572\endref

\ref
\key KL\by \ B.~Kleiner and J.~Lott\paper Notes on Perelman's papers,
http://arxiv.org/abs/math/0605667\endref

\ref
\key KZ\by \ S.~Kuang and Q.S.~Zhang\paper A gradient estimate for
all positive solutions of the conjugate heat equation under Ricci
flow, http://arxiv.org/abs/math/0611298\endref

\ref
\key LSU\by \ \ O.A.~Ladyzenskaya, V.A.~Solonnikov, and
N.N.~Uraltceva\book Linear and quasilinear equations of
parabolic type\publ Transl. Math. Mono. Vol 23,
Amer. Math. Soc.\publaddr Providence, R.I.\yr 1968\endref

\ref
\key LK\by \ P.~Li and L.~Karp\paper The heat equation on complete 
Riemannian manifolds (unpublished manuscript)\endref

\ref
\key MT\by \ \ J.W.~Morgan and G.~Tian\paper Ricci flow and the Poincar\'e
Conjecture, http://arXiv.org/abs/math.DG\linebreak /0607607\endref

\ref
\key NT\by \ L.~Ni and L.F.~Tam\paper K\"ahler-Ricci flow and the 
Poincar\'e-Lelong equation\jour Comm. Anal. and Geom.\vol 12(1)\yr 2004
\pages 111-141\endref

\ref
\key P\by G.~Perelman\paper The entropy formula for the Ricci flow and its 
geometric applications,\linebreak http://arXiv.org/abs/math.DG/0211159\endref 

\ref
\key S1\by W.X.~Shi\paper Deforming the metric on complete Riemannian 
manifolds\jour J. Differential Geom.\vol 30\yr 1989\pages 223--301\endref

\ref
\key S2\by W.X.~Shi\paper Ricci deformation of the metric on complete 
non-compact Riemannian manifolds\jour J. Differential Geom.\vol 30\yr 1989
\pages 303--394\endref

\ref
\key S3\by W.X.~Shi\paper Ricci flow and the uniformization on complete 
noncompact K\"ahler manifolds\jour J. Differential Geom.\vol 45\yr 1997
\pages 94--220\endref

\ref
\key W\by J.~Wang\paper The heat flow and harmonic maps between complete
manifolds\jour J. Geometric Analysis\vol 8(3)\yr 1998\pages 485--514
\endref

\endRefs
\enddocument